\input epsfx.tex
\input amstex
\catcode`@=12

\frenchspacing
\documentstyle{amsppt}
\magnification=\magstep1
\baselineskip=14pt
\vsize=18.5cm
\footline{\hfill\sevenrm version 20090202}
\def\C{\bold C}
\def\H{\bold H}
\def\R{\bold R}
\def\Q{\bold Q}
\def\Z{\bold Z}
\def\F{\bold F}
\def\Fp{{\F}_p}

\def\P{{\bold P}}
\def\O{{\Cal O}}
\def\Pic{\text{\rm Pic}}
\def\Mat{\text{\rm Mat}}
\def\Jac{\text{\rm Jac}}
\def\Sp{\text{\rm Sp}}
\def\GL{\text{\rm GL}}
\def\SL{\text{\rm SL}}
\def\Im{\text{\rm Im}}
\def\Tr{\text{\rm Tr}}
\def\Spec{\text{\rm Spec}}
\def\Am{{\Cal A}} 
\def\G2#1{{\Gamma^{(2)}#1}}
\def\Gs2#1{{\Gamma_0^{(2)}#1}}
\def\Ys2#1{{Y_0^{(2)}#1}}

\def\endproof{\hfill$\square$\par}
\def\isar{\ \smash{\mathop{\longrightarrow}\limits^{\thicksim}}\ }
\gdef\mapright#1{\ \smash{\mathop{\longrightarrow}\limits^{#1}}\ }
\def\eqdef{{\thinspace\displaystyle\mathrel{\mathop=^{\text{\rm def}}}\thinspace}}
\newcount\refCount
\def\newref#1 {\advance\refCount by 1
\expandafter\edef\csname#1\endcsname{\the\refCount}}
\newref AZ    
\newref Art   
\newref BB    
\newref Var   
\newref Lange 
\newref BM    
\newref CN    
\newref CrMan 
\newref GK    
\newref Du    
\newref EZ    
\newref EL    
\newref Mor   
\newref GT    
\newref PGH   
\newref pAdic 
\newref GS    
\newref Igu   
\newref Webb  
\newref Koe   
\newref ShaB   
\newref Mes   
\newref Mura  
\newref Schb  
\newref Scha  
\newref Sha   
\newref ShaA   
\newref Suth  

\topmatter
\title
Modular polynomials for genus 2
\endtitle
\author Reinier Br\"oker, Kristin Lauter
\endauthor
\address
Microsoft Research, One Microsoft Way, Redmond, WA 98052, USA
\email reinierb\@microsoft.com, klauter\@microsoft.com \endemail
\endaddress
\abstract
Modular polynomials are an important tool in many algorithms involving elliptic
curves. In this article we investigate their generalization to the genus 2
case following pioneering work by Gaudry and Dupont. We prove various
properties of these genus 2 modular polynomials and give an improved way
to explicitly compute them.
\endabstract
\endtopmatter

\document

\head 1. Introduction
\endhead
\noindent
The `classical' modular polynomial $\Phi_N \in \Z[X,Y]$ was introduced
by Kronecker more than 100 years ago. The polynomial $\Phi_N$ is a model
for the modular curve $Y_0(N)$ parametrizing cyclic $N$-isogenies
between elliptic curves. As is shown in the examples below, the explicit
computation of $\Phi_N$ has led to various speed ups in algorithms
using elliptic curves.\par
\ \par\noindent
{\bf Examples.}
1. Schoof's original algorithm [\Scha] to count the number of points on
an elliptic curve $E/\Fp$
was rather impractical as one had to compute with the {\it complete\/}
$l$-torsion of $E$ for various small primes $l \not = p$. The key to the
improvements made by Atkin and Elkies [\Schb, Sections 6--8] is to work with
a one-dimensional eigenspace
of $E[l] \cong \Z/l\Z \times \Z/l\Z$. Instead of using division polynomials
of degree $(l^2-1)/2$ one can now use the modular polynomial $\Phi_l$ of
degree $l+1$. The `Schoof-Elkies-Atkin'-algorithm behaves very well in
practice, and primes $p$ of several thousand digits are now feasible [\Mor].

2. Both the primality proving algorithm ECPP [\GK, Section 14D] and
efficient constructions of cryptographically secure elliptic curves
[\CrMan, Chapter 18] rely on the computation
of the Hilbert class polynomial $H_\O$ for a certain imaginary quadratic
order $\O$. At the moment, the fastest known algorithm to compute $H_\O
\in \Z[X]$ is to first compute its reduction modulo various primes~$p$ and then
apply the Chinese remainder theorem [\Var, \Suth].
If we know one root of $H_\O$ in $\Fp$,
then we can compute the other roots by exploiting the modular polynomials
$\Phi_l$ for small primes $l$ generating the class group of $\O$. This
observation is crucial to the practical behaviour of the `CRT-algorithm'.\par
\ \par\noindent
The algorithms in the examples above have analogues for genus 2, see
[\PGH, \EL]. Whereas
the situation in well understood in genus 1 and algorithms are fast, the
computations are still in their infancies in genus 2. Typically, algorithms
only terminate in a reasonable amount of time for very small examples.
Inspired by the speed ups modular polynomials give in the genus 1 case, we
investigate modular polynomials for genus 2 in this paper.

In [\GS], Gaudry and Schost examine a tailor-made variant of $\Phi_N$ in
genus 2 to improve point counting on genus 2 curves over finite fields. The polynomial they construct has
factorization properties similar to $\Phi_N$, and this enables them to speed
up point counting in genus 2 in the spirit of the improvements made by Atkin
and Elkies in the genus 1 case.

In this paper, we consider a `direct' generalization of $\Phi_N$ to
genus~2. Such a generalization was first defined by Gaudry in his
PhD-thesis [\GT, Ch.\ 3]. Whereas we have
one polynomial $\Phi_N$ in the genus 1 case, we now get $3$ polynomials
$P_N, Q_N, R_N$ for every~$N\geq2$.
At the time, it was not possible
to {\it explicitly compute\/} these polynomials $P_N, Q_N, R_N$ in
the simplest case $N=2$. Dupont considered the problem of computing the
modular polynomials in this thesis [\Du, Sec.\ 10.4] as well. He was able to compute
them for $N=2$, and he gives some partial results for $N=3$. The resulting
polynomials are {\it huge\/}. 
Nevertheless, knowing just a few modular polynomials
will speed up the `CRT-algorithm' [\EL] to compute class fields
of degree~4 CM-fields. This in turn will lead to faster algorithms to
construct cryptographically secure Jacobians of hyperelliptic curves.
This paper solely focuses on the definitions and properties of modular
polynomials for genus 2 however.

We reconsider the polynomials defined by Gaudry and Dupont in this article. We
will mostly restrict to the case that $N=p$ is prime.
By combining the function field approach of Gaudry with finite symplectic
geometry over $\Fp$, we are able to prove Gaudry's conjecture on the degree of the
modular polynomials. 
Furthermore, we prove that the polynomials $P_p, Q_p,
R_p$ have coefficients in $\Q(j_1,j_2,j_3)$. We also give a heuristic in 
Section~6 implying that the bit size of the polynomial $P_p$ grows 
like~$p^{12}$. We extend Dupont's results in
the following way: the algorithm in [\Du] requires `guessing' the
denominator of the coefficients of $P_p,Q_p,R_p$. In Section~6 we prove that
the denominator is closely linked to the {\it Humbert surface\/} of
discriminant~$p^2$ describing
$(p,p)$-split Jacobians. Our result considerably helps the computation.

The structure of this article is as follows.
After recalling some classical facts of polarized abelian
varieties in Section 2, we consider level structures in Section~3, and
for primes $p$ we define the $2$-dimensional analogue $\Ys2(p)$ of the curve
$Y_0(p)$. In Section~4 we consider functions on $\Ys2(p)$ and this leads
naturally to the definition of the genus~2 modular polynomials.
 We explain the moduli interpretation behind $P_p, Q_p, R_p$, and compute
their degrees. By studying the Fourier coefficients of the Igusa $j$-functions,
 we show in Section~5 that the genus 2 modular polynomials have
rational coefficients. Section 6 focuses on the {\it computation\/}
of the modular polynomials. In the special case $p=2$, the modular polynomials
are closely linked to the {\it Richelot isogeny\/}. We explain the precise
relation in Section 7.

\head 2. Polarized abelian varieties
\endhead
\noindent
We recall some classical facts about complex abelian varieties. For modular
polynomials we are only interested in the 2-dimensional case, but as much
of the theory generalizes, we work in arbitrary dimension $g\geq 1$ in this
section.

Let $A/\C$ be a $g$-dimension abelian variety, and write $A^\vee = \Pic^0(A)$
for its dual. It is well known that every abelian variety admits a
{\it polarization\/}, i.e., an isogeny $\varphi: A \rightarrow A^\vee$. If
$\varphi$ is an isomorphism, we call $A$ principally polarized. The
Jacobian of a curve is the classical example of a principally polarized
abelian variety: the choice of a base point on the curve determines the
polarization.

The complex points $A(\C)$ have the natural structure of a $g$-dimensional
complex torus $\C^g/L$. Here, $L \subset \C^g$ is a full lattice. In order
to characterize those complex tori $\C^g/L$ that arise as abelian varieties, we
define the polarization of a torus. We follow the classical approach via
Riemann forms. A skew-symmetric form $E: L \times L \rightarrow \Z$ can be
extended to a form $E: \C^g \times \C^g \rightarrow \R$, and we call $E$
a {\it Riemann form\/} if
\medskip
\item{{\bf 1.}} $E(x,y) = E(ix,iy)$ for all $x,y \in \C^g$
\item{{\bf 2.}} the Hermitian form $H(x,y) = E(ix,y) + iE(x,y)$ is positive
definite. \par\medskip\noindent
A torus $\C^g/L$ is called polarizible if it admits a Riemann form. The link
with the definition from the previous paragraph is that the
Hermitian form $H$ from condition 2 defines a map
$$
\varphi_E: (\C^g/L) \mapright{} (\C^g/L)^\vee
$$
sending $x$ to $H(x,\cdot)$. If $\varphi_E$ is an isomorphism, we say that the
torus $\C^g/L$ is principally polarized.
The following theorem describes the precise connection between complex abelian
varieties and polarizable tori.\par
\noindent
\proclaim{Theorem 2.1}{The category of complex abelian varieties is equivalent
to the category of polarizable tori via the functor $A \rightarrow A(\C)$.}
\endproclaim
\noindent
{\bf Proof.} See [\Lange, Chapter~4]. \endproof\par
\ \par\noindent
Let $L$ be a full lattice in $\C^g$ such that the torus $\C^g/L$ is
principally polarizable. One can choose a basis of $L$ such that $L$ is given by
$\Z^g + \Z^g\tau$ with $\tau$ a complex $g\times g$-matrix. The fact that
the lattice admits a Riemann form now
translates into the fact that $\tau$ is an element of the $g$-dimensional
{\it Siegel upper half plane\/}
$$
\H_g = \{ \tau \in \Mat_{g}(\C) \mid \tau^T = \tau, \Im(\tau) \text{\rm\ positive definite} \}.
$$
Conversely, for $\tau \in \H_g$ the torus $\C^g/(\Z^g+\Z^g\tau)$ is
principally polarizable. The Hermitian form given by $H(x,y) = x \Im(\tau)^{-1}
\overline y^T$ is by construction positive definite, and the Riemann form is
given by $E(x,y) = \Im(H(x,y))$.

Next we consider what happens if we change a basis for the $2g$-dimensional
polarized lattice~$L$. By the `elementary divisor theorem', we can choose a
basis for $L$ such that the Hermitian form $H: L \times L \rightarrow \Z$ is
given by the matrix
$$
J = \Bigl( {0 \atop -1_g} \thinspace {1_g \atop 0} \Bigr)
$$
where $1_g$ denotes the $g\times g$ identity matrix.
The general linear group $\GL(2g,\Z)$ stabilizes the lattice $L$, and
the subgroup
$$
\Sp(2g,\Z) = \{ M \in \GL(2g,\Z) \mid M J M^T = J \} \subseteq \GL(2g,\Z)
$$
that respects the Hermitian form is called the {\it symplectic group\/}.
Explicitly, a matrix $M = \Bigl( {a \atop c} \thinspace {b \atop d}
\Bigr)\in \GL(2g,\Z)$ is symplectic if and only if the $g\times g$
matrices $a,b,c,d$ satisfy the relations $ab^T = b^Ta$, $cd^T = d^Tc$
and $ad^T-bc^T = 1_g$.

For $\Bigl( {a \atop c} \thinspace {b \atop d} \Bigr) \in \Sp(2g,\Z)$ and
$\tau \in \H_g$, the $g\times g$-matrix $c\tau +d$ is invertible. Indeed,
if $c\tau +d$ had determinant zero, there would be an element
$y\in \C^g$ with $y\Im(\tau) y^T=0$, contradicting that $\Im(\tau)$ is
positive definite. We define an action of the symplectic group $\Sp(2g,\Z)$ on
the Siegel upper half plane $\H_g$ by
$$
\Bigl( {a \atop c} \thinspace {b \atop d} \Bigr)\tau = {a\tau + b \over c\tau+d}, \eqno(2.1)
$$
where dividing by $c\tau +d$ means multiplying on the right with the multiplicative inverse
of the $g\times g$-matrix $c\tau +d$. An explicit check shows that the
right hand side of (2.1) indeed lies in $\H_g$.

The map
$$
\tau \mapsto A_\tau = \C^g/(\Z^g+\Z^g\tau).
$$
induces a canonical bijection between the quotient space
$\Am_g = \Sp(2g,\Z) \backslash \H_g$ and the set of isomorphism classes of
principally polarized $g$-dimensional abelian varieties. In fact, the
space $\Am_g$ is a {\it coarse\/} moduli space for principally polarized
abelian varieties of dimension~$g$.

We close this section by zooming in on the 1-dimensional case, i.e., the
case of elliptic curves. An elliptic curve is isomorphic to its dual, so
polarizations do not play a real role here. Indeed, every complex torus
is polarizable, and the Siegel space $\H_1$ equals the Poincar\'e upper
half plane~$\H$. The symplectic group $\Sp(2,\Z)$ equals the special
linear group $\SL_2(\Z)$ in this case. The space $\Sp_2(\Z) \backslash \H$ is
in canonical bijection with the set of isomorphism classes of elliptic
curves.
\head 3. Isogenies
\endhead
\noindent
Let $A/\C$ be a 2-dimensional principally polarized abelian variety, and let
$N \geq 1$ be a positive integer. The $N$-torsion $A[N]$ of $A$ is,
non-canonically,
isomorphic to $(\Z/N\Z)^4$. The polarization on $A$ induces a symplectic
form $v$ on the rank~4 $(\Z/N\Z)$-module $A[N]$.
We choose a basis for $A[N]$ such that $v$ is given by the matrix
$$
\Bigl( {0 \atop -1_2} \thinspace {1_2 \atop 0} \Bigr),
$$
and we let
$\Sp(4,\Z/N\Z)$ be the subgroup of the matrix group $\GL(4,\Z/N\Z)$ that
respects~$v$. A subspace $G \subset A[N]$ is called {\it isotropic\/} if
$v$ restricts to the zero-form on $G \times G$, and we say that
$A$ and $A'$ are $(N,N)$-isogenous if there is an isogeny $A \rightarrow A'$
whose kernel is isotropic of order $N^2$.

The full congruence subgroup $\G2(N)$ of level $N$ is defined as the kernel
of the reduction map $\Sp(4,\Z) \rightarrow \Sp(4,\Z/N\Z)$. Explicitly,
a matrix $\bigl( {a \atop c} \thinspace {b \atop d} \bigr)$ is contained
in $\G2(N)$ if and only if we have $a,d \equiv 1_2 \bmod N$ and $b,c
\equiv 0_2 \bmod N$. The congruence subgroup $\G2(N)$ fits in an exact sequence
$$
1 \mapright{} \G2(N) \mapright{} \Sp(4,\Z) \mapright{} \Sp(4,\Z/N\Z)
\mapright{} 1.
$$
The surjectivity is not completely trivial, see [\AZ, Lemma 3.2].

The 2-dimensional analogue of the subgroup $\Gamma_0(N) \subset \SL_2(\Z)$
occuring in the equality $Y_0(N) = \Gamma_0(N) \backslash \H$ of Riemann
surfaces is the group
$$
\Gs2(N) = \left\{ \Bigl( {a \atop c} \thinspace {b \atop d} \Bigr) \in
\Sp(4,\Z) \mid c \equiv 0_2 \bmod N \right\}.
$$
From now on, we restrict to the case $N=p$ prime. The reason for this
restriction is that the finite symplectic geometry we need in the remainder
of this section is much easier for vector spaces over finite fields than
for modules over arbitrary finite rings.
The following lemma gives the link between the group $\Gs2(p)$ and isotropic
subspaces of the $p$-torsion.

\proclaim{Lemma 3.1}
The index $[\Sp(4,\Z) : \Gs2(p)]$ equals the number of $2$-dimensional
isotropic subspaces of the $\Fp$-vector space $\Fp^4$.
\endproclaim
\noindent{\bf Proof.} We map $\Gs2(p)$ to a subgroup $H \subset \Sp(4,\Fp)$.
The inclusion $\G2(p) \subseteq \Gs2(p)$ shows that we have $[\Sp(4,\Z):\Gs2(p)]
= [\Sp(4,\Fp):H]$. The group $H$ is parabolic, and occurs
as stabilizer of the 2-dimensional isotropic subspace
$\Fp \times \Fp \times 0 \times 0$ of the symplectic space $\Fp^4$.
The group $\Sp(4,\Fp)$ permutes the $2$-dimensional
isotropic subspaces transitively by Witt's extension theorem
[\Art, Theorem 3.9], and the lemma follows. \endproof
\ \par\noindent
Let $S(p)$ be the set of equivalence classes of pairs $(A,G)$, with $A$ a
2-dimensional principally polarized abelian variety and $G \subset A[p]$ a
2-dimensional isotropic subspace. Here, two pairs $(A,G)$ and $(A',G')$
are said to be isomorphic if there exists an isomorphism of abelian varieties
$\varphi: A \rightarrow A'$ with $\varphi(G) = G'$.

\proclaim{Theorem 3.2}
The quotient space $\Gs2(p) \backslash \H_2$ is in canonical bijection with
the set $S(p)$ via $\Gs2(p)\tau \mapsto (A_\tau, \langle (1/p,0,0,0),
(0,1/p,0,0) \rangle)$ where $A_\tau = \C^2/(\Z^2+\Z^2\tau)$ is the variety
associated to $\tau$.
\endproclaim
\noindent{\bf Proof.} The group $\Gs2(p)$ stabilizes the subspace
$G = \langle (1/p,0,0,0), (0,1/p,0,0)\rangle$ of $A[p]$, and the image
$(A_\tau,G)$ therefore does not depend on the choice of a representative $\tau$.

If $A_\tau$ and $A_{\tau'}$ are isomorphic, then there exists
$\psi\in\Sp(4,\Z)$ with $\psi\tau = \tau'$. If an isomorphism $A_\tau \isar
A_{\tau'}$ maps the group $G$ to $G' = \langle (1/p,0,0,0),
(0,1/p,0,0)\rangle \subset A_{\tau'}[p]$, then $\psi$ lies in $\Gs2(p)$. Hence,
our map is injective.

To prove surjectivity, we first note that every 2-dimensional principally
polarized abelian variety occurs as some $A_\tau$ by Theorem 2.1. The
theorem now follows directly from Lemma 3.1.\endproof
\ \par\noindent
As a quotient space, the 2-dimensional analogue of the curve $Y_0(p)$ is
$$
\Ys2(p) \eqdef \Gs2(p) \backslash \H_2.
$$
We close this section by showing how to give $\Ys2(p)$ the structure of a
quasi-projective variety. Siegel defined a metric on $\H_2$, a
generalization of the Poincar\'e metric in dimension 1, that
respects the action of the symplectic group. With this metric, $\Ys2(p)$
becomes a topological space. Just as in the 1-dimensional case $Y_0(p)$, it
is not compact. There are several ways to compactify it, one of which is
the {\it Satake compactification\/}
$$
\Ys2(p)^* = \Ys2(p) \cup Y_0(p) \cup \P^1(\Q).
$$
By the Baily-Borel theorem [\BB], the space $\Ys2(p)^*$ has a natural
structure as a projective variety $V$, and the space $\Ys2(p)\subset V$ is
therefore naturally a quasi-projective variety.
\head 4. Functions
\endhead
\noindent
The left-action of $\Sp(4,\Z)$ on the Siegel upper half plane $\H_2$ induces a
natural right-action on the set of functions from $\H_2$ to $\P^1(\C)$ via
$(fM)(\tau) = f(M\tau)$. The fixed points under this action are called
{\it rational Siegel modular functions\/}. Igusa defined [\Igu, Theorem~3]
three algebraically independent rational Siegel modular functions
$\H_2 \rightarrow \P^1(\C)$
that generate the function field $K$ of $\Am_2 = \Ys2(1)$. The functions
$j_1,j_2,j_3$ that most people use nowadays are slightly different from
Igusa's and we recall their definition first.

Let $E_k(\tau) = \sum_{(c,d)}(c\tau+d)^{-k}$ be the 2-dimensional Eisenstein
series. Here, the sum ranges over all co-prime symmetric $2\times 2$-integer
matrices that are non-associated with respect to left-multiplication
by $\GL(2,\Z)$. We define $$\chi_{10} = -43867\cdot 2^{-12} \cdot 3^{-5} \cdot
5^{-2} \cdot 7^{-1} \cdot 53^{-1} (E_4 E_6 -E_{10})$$ and $$\chi_{12} =
131\cdot 593\cdot 2^{-13} \cdot 3^{-7} \cdot 5^{-3} \cdot 7^{-2} \cdot 337^{-1}
(3^2\cdot 7^2 E_4^3 +2 \cdot 5^3 E_6^2 -691 E_{12}),$$ where the constants
in $\chi_{10}$ and $\chi_{12}$ should be regarded as `normalization factors'.
We then define
$$
j_1 = 2\cdot 3^5 {\chi_{12}^5 \over \chi_{10}^6}, \qquad j_2 = 2^{-3} 3^3
{E_4 \chi_{12}^3 \over \chi_{10}^4}, \qquad j_3 = 2^{-5} \cdot 3 {E_6
\chi_{12}^2 \over \chi_{10}^3} + 2^{-3}\cdot 3^2 {E_4 \chi_{12}^3 \over \chi_{10}^4}.$$
We have $K = \C(j_1,j_2,j_3)$ and via the moduli
interpretation for $\Am_2$, the Igusa
functions are functions on the set of principally polarized 2-dimensional
abelian varieties.  The functions $j_1,j_2,j_3$ have poles at $\tau
\in \H_2$ corresponding to products of elliptic curves with the product
polarization.

For a fixed prime $p$, we define three functions
$$
j_{i,p} : \H_2 \rightarrow \P^1(\C) \eqno(4.1)
$$
by $j_{i,p}(\tau) = j_i(p\tau)$. These functions arise naturally in the
study of $(p,p)$-isogenous abelian varieties as we have
$$
j_{i,p}(A_\tau / \langle (1/p,0,0,0), (0,1/p,0,0) \rangle) = j_i(A_{p\tau})
= j_{i,p}(\tau).
$$
\proclaim{Lemma 4.1}
The functions $j_{i,p}$ defined in (4.1) above are invariant under the action
of~$\Gs2(p)$.
\endproclaim
\noindent
{\bf Proof.}
Let $M = \bigl( {a \atop c} \thinspace {b \atop d} \bigr)$ be
an element of $\Gs2(p)$, and write $c = pc'$ with $c'$ a $2\times 2$-matrix with
integer coefficients. We compute
$$
j_{i,p}(M\cdot\tau) = j_i(p(M\cdot\tau)) = j_i((pa\tau + pb)/(c\tau +d)) = j_i(B
\cdot p\tau),
$$
with $B = \bigl( {a \atop c'} \thinspace {pb \atop d} \bigr)$. It is clear
that $B$ is again symplectic, so we have $j_i(B \cdot p\tau) = j_{i,p}(\tau)$.
\endproof\par
\ \par\noindent
The functions $j_{i,p}$ have poles at $\tau\in\H_2$ corresponding to
$(p,p)$-{\it split\/} principally polarized abelian varieties, i.e., varieties
that are $(p,p)$-isogenous to a product of elliptic curves with the product
polarization.\par
\proclaim{Lemma 4.2}
For a prime $p$, the function field of $\Ys2(p)/\C$ equals $K(j_{i,p})$ for
every $i=1,2,3$.
\endproclaim
\noindent
{\bf Proof.} The function field of $\Am_2=\Ys2(1)$ equals $K=\C(j_1,j_2,j_3)$
and the function field $\C(\Ys2(p))$ is an extension of $K$ of degree
$[\Sp(4,\Z):\Gs2(p)]$. The functions $j_{i,p}$ are contained in
$\C(\Ys2(p))$ by Lemma 4.1. It suffices to show that, for fixed~$i$, the
functions $\{j_{i,p}(\alpha\tau)\}_{\alpha\in\Gs2(p)\backslash\Sp(4,\Z)}$ are
distinct. If two of these functions are equal, then the stabilizer
$S \subset \Sp(4,\Z)$ of $j_{i,p}$ inside $\Sp(4,\Z)$ {\it strictly\/}
contains $\Gs2(p)$. The images of $S$ and $\Gs2(p)$ under the reduction
map $\pi: \Sp(4,\Z) \twoheadrightarrow \Sp(4,\Fp)$ then satisfy
$$
\pi(S) \supsetneq \pi(\Gs2(p)).
$$
The group $\pi(\Gs2(p))$ is the stabilizer of an isotropic subspace of
$\Sp(4,\Fp)$ and is therefore {\it maximal\/} by [\Webb, Theorem~4.2].
Hence, $\pi(S)$ equals the
full group $\Sp(4,\Fp)$ and $S$ has to equal $\Sp(4,\Z)$. This is absurd.
\endproof\par
\ \par\noindent
The $p$th {\it modular polynomial\/} $P_p$ for $j_1$ is defined as the
minimal polynomial of $j_{1,p}$ over $K=\C(j_1,j_2,j_3)$. It has
degree $[\Sp(4,\Z) : \Gs2(p)]$ and its coefficients are rational functions
in $j_1,j_2,j_3$ with complex coefficients. The evaluation map
$\varphi_\tau: \C(j_1,j_2,j_3) \rightarrow \C$ sending $j_i$ to $j_i(\tau)$
maps $P_p$ to a polynomial $P_{p,\tau} \in \C[X]$. The roots of
$P_{p,\tau}$ are the $j_1$-invariants of principally polarized abelian
varieties that are $(p,p)$-isogenous to a variety with
$j$-invariants $j_1(\tau), j_2(\tau), j_3(\tau)$.

The functions $j_{2,p},j_{3,p}$ are contained in $K(j_{1,p}) = K[j_{1,p}]$
and we define $R_p, Q_p\in \C(j_1,j_2,j_3)[X]$ to be the monic
polynomials of degree less than $\deg(P_p)$ satisfying
$$
j_{2,p} = R_p(j_{1,p}) \qquad\qquad
j_{3,p} = Q_p(j_{1,p}). \eqno(4.2)
$$
The evaluation map $\varphi_\tau$ maps $Q_p,R_p$ to
polynomials $Q_{p,\tau},R_{p,\tau} \in \C[j_{1,p}]$. If $x\in\C$ is a root of
$P_{p,\tau}$, then
$$
(x,Q_{p,\tau}(x),R_{p,\tau}(x))
$$
are $j$-invariants of a principally polarized abelian variety that is
$(p,p)$-isogenous to a variety with invariants $j_1(\tau),j_2(\tau),j_3(\tau)$.

\head 5. Field of definition
\endhead
\noindent
A holomorphic map $\psi: \H_2 \rightarrow \C$ is called a {\it Siegel modular
form\/} of degree $w \geq 0$ if it satisfies the functional equation
$$
\psi(\left( {a \atop c} \thinspace {b \atop d} \right) \tau) = \det(c\tau+d)^{w} \psi(\tau)
$$
for all matrices $\bigl( {a \atop c} \thinspace {b \atop d} \bigr) \in
\Sp(4,\Z)$. The Eisenstein series $E_w$ are Siegel modular forms of degree~$w$.
Any Siegel modular form is invariant under the transformation $\tau \mapsto
\tau+b$ and therefore admits a {\it Fourier expansion\/}
$$
\psi = \sum_{T}a(T) \exp(2\pi i \Tr(T\tau)),
$$
where the summation ranges over all $2\times 2$ symmetric `half-integer'
matrices, i.e., symmetric matrices with integer entries on the diagonal and
off-diagonal entries in~${1 \over 2}\Z$. The coefficients $a(T)$ are called the
Fourier coefficients of~$\psi$. As discovered by Koecher [\Koe], they are zero
unless $T$ is positive semi-definite.
For $T = \bigl( {a \atop b/2} \thinspace {b/2 \atop c} \bigr)$ and $\tau =
\bigl( {\tau_1 \atop \tau_2} \thinspace {\tau_2 \atop \tau_3} \bigr)$ we have
$$
\Tr(T\tau) = a\tau_1 + b\tau_2 + c\tau_3.
$$
Writing $q_i = \exp(2\pi i \tau_i)$, we see that we can express a modular
form as $\psi = \sum_{k,l,m} c_{k,l,m} q_1^k q_2^l q_3^m$. By Koecher's
result, the summation ranges over non-negative $k,l$ and $m$ satisfying
$4m-kl\geq 0$.

A Siegel modular form is called a {\it cusp form\/} if the Fourier coefficients
$a(T)$ are zero for all $T$ that are semi-definite but not definite. One of
the classical examples of a cusp form is
$$
\chi_{10} = -43867\cdot 2^{-12} \cdot 3^{-5} \cdot 5^{-2} \cdot 7^{-1}
\cdot 53^{-1} (E_4 E_6 -E_{10}),
$$
which appears in the denominator of the Igusa $j$-functions. If we express
$\chi_{10}$ in its `$q_i$-expansion', every term is divisible by
$q_1 q_2 q_3$. The `normalization factor' ensures that the $q_1 q_2 q_3$-term
has coefficient~1.\par
\ \par
\proclaim{Lemma 5.1} The Igusa functions $j_i$ have a Laurent series expansion
in $q_1$, $q_2$, $q_3$ with rational coefficients.
\endproclaim
\noindent
{\bf Proof.} The denominator of all three Igusa functions is a constant
multiple of a power of the cusp form $\chi_{10}$. The product $(q_1q_2q_3)^{-1} \chi_{10}$
has a {\it non-zero\/} constant term and is therefore invertible in the ring
$\C[[q_1,q_2,q_3]]$. This shows that the Igusa functions have a Laurent
series expansion.

As the Fourier coefficients of the Eisenstein series are rational [\EZ, Corollary~2 to Theorem~6.3], the
coefficients of the Laurent expansion of $j_i$ are rational. \endproof\par
\ \par\noindent
In genus 1, it is not hard to prove that the Fourier coefficients of the
$j$-function are rational. A deeper result is that they are {\it integral\/}.
This is no longer true in genus~2: the coefficients of the expansion of the
Igusa functions have `true' denominators.\par
\ \par
\noindent
\proclaim{Theorem 5.2} For any prime~$p$, the modular polynomials
$P_p,Q_p,R_p$ lie in the ring $\Q(j_1,j_2,j_3)[X]$.
\endproclaim
\noindent
{\bf Proof.} We only give the proof for $P_p$, the proof for $Q_p$ and $R_p$ is
highly similar. We can write
$$
P_p = \sum_{m\geq 0} {\sum_{a,b,c} c_{m,a,b,c} j_1^a j_2^b j_3^c \over
                      \sum_{a,b,c} d_{m,a,b,c} j_1^a j_2^b j_3^c} X^m, \eqno(5.1)
$$
and we have to prove that the coefficients $c_{m,a,b,c}$ and $d_{m,a,b,c}$ are
rational. We substitute the Laurent series expansion of $j_1,j_2,j_3,j_{1,p}$
into the equation $P_p = 0$. By equating powers of $q_1^a q_2^b q_3^c$, we
get a set of linear equations for the $c_{m,a,b,c}$ and $d_{m,a,b,c}$.

Over the complex numbers, this system of equations has a unique solution. As
the coefficients of the equations are rational by Lemma 5.1, this solution
must be rational.\endproof
\ \par\noindent
{\bf Remark 5.3.} We can reduce the polynomials $P_p, Q_p, R_p$ modulo a
prime $l$. A natural question is if these reduced polynomials still satisfy
a moduli interpretation as in (4.2). Whereas reduction of modular curves is
relatively well understood, the situation is more complicated for general
Siegel modular varieties. In our situation, the answer is given by a theorem
of Chai and Norman [\CN, Corollary 6.1.1]. They look at the algebraic stack
${\Cal A}_{2,\Gs2(p)}$ and prove that the structural morphism
${\Cal A}_{2,\Gs2(p)} \rightarrow \Spec\ \Z$ is faithfully flat,
Cohen-Macaulay and {\it smooth outside~$p$.} Concretely, this means that
the moduli interpretation (4.2) remains valid modulo primes~$l \not = p$.

\head 6. Explicit computations
\endhead
\noindent
In this Section we give a method to compute the modular polynomials
$P_p, Q_p, R_p \in \Q(j_1,j_2,j_3)$ and indicate what the computational
difficulties are. We begin with the degree of $P_p$.\par
\ \par
\proclaim{Lemma 6.1} For a prime $p$, we have $[\Sp(4,\Z):\Gs2(p)] =
(p^4-1)/(p-1)$.
\endproclaim
\noindent
{\bf Proof.} By Lemma 3.1, we have to count the number of 2-dimensional
isotropic subspaces of the symplectic space $V=\Fp^4$.

Any two-dimensional isotropic subspace of $V$ contains $(p^2-1)/(p-1) =
p+1$ lines. Conversely, a line $l \subset V$ is contained in $(p+1)$ isotropic
subspaces. To see this, we note that we can need to select a second line $m$
such that $\langle l, m \rangle$ is isotropic of dimension~$2$. The complement
of $l$ is 3-dimensional, and out of the $(p^3-1)/(p-1)$ lines, only
$(p^2-1)/(p-1) = p+1$ yield an {\it isotropic\/} subspace.

We see that the number of 2-dimensional isotropic subspaces equals the number
of lines in $V$. This yields the lemma.\endproof\par
\ \par\noindent
{\bf Remark.} It is easy to give coset representatives for $S =
\Sp(4,\Z)/\Gs2(p)$. In genus~1, we can take the set
$$
\left\{ \left( {1 \atop i} \thinspace {0\atop 1} \right) \mid i \in \Fp
\right\} \cup \left\{ \left( {0 \atop 1} \thinspace {-1 \atop 0} \right)
\right\} \eqno(6.1)
$$
as a set of coset representatives for $\SL(2,\Z)/\Gamma_0(p)$. Inspired by the
set in (6.1) we write down
$$
\left\{ \left( {1_2 \atop \bigl({a \atop b} \thinspace {b \atop c}\bigr)}
\thinspace {0_2 \atop 1_2}\right) \mid a,b,c \in\Fp\right\} \cup \left\{
\left( {0_2 \atop 1_2} \thinspace {-1_2 \atop \bigl( {a \atop b} \thinspace
{b \atop c} \bigr)} \right) \mid  ac = b^2 \in \Fp \right\}, \eqno(6.2)
$$
a set of cardinality $p^3+p^2$. We are missing $p+1$ matrices. One can check
that
$$
\left\{ \left( { {1 \atop 0} \thinspace {0 \atop 0} \atop {0 \atop -a}
\thinspace {0 \atop 1} } \thinspace { { 0 \atop 0} \thinspace {0 \atop -1} \atop
{1 \atop 0} \thinspace {a \atop 0} } \right) \mid a \in \Fp \right\} \cup
\left\{ \left( { {-1 \atop 0} \thinspace {-1 \atop 0} \atop  {0 \atop 1}
\thinspace {0 \atop 0} } \thinspace { {1 \atop -1} \thinspace {-1 \atop 1}
\atop {0 \atop 0} \thinspace {-1 \atop -1} }\right) \right\} \eqno (6.3)
$$
is a set of $p+1$ matrices that is independent of the set in (6.2). We
note that this is the same set of matrices that Dupont found in his
thesis [\Du].\par
\ \par
\noindent
To compute $R_p$ and $Q_p$, we have to write $j_{2,p}$ and $j_{3,p}$ as
rational functions in $j_1,j_2,j_3$ and $j_{1,p}$. For $M$ ranging over the
cosets $S = \Sp(4,\Z) / \Gs2(p)$, the functions $j_{i,p}(M\tau)$ are distinct by
Lemma~4.2. Inspired by the formulas in [\pAdic, Section 7.1] we note that, by
Lagrange interpolation, the polynomials
$$
F_{k,p}(X) = \sum_{M\in S} \Biggl( \prod_{B\in S \atop B\not = M} {X -
    j_{1,p}(B\tau) \over j_{1,p}(M\tau) - j_{1,p}(B\tau)}\Biggr) j_{2,p}(M\tau)
$$
satisfy $F_{k,p}(j_{1,p}(C\tau)) = j_{2,p}(C\tau)$ for $k=2,3$ and all
$C\in S$. As the coefficients of $F_{k,p}$ are, by construction, invariant
under the action of $\Sp(4,\Z)$ the polynomials $F_{k,p}$ are contained in
$\Q(j_1,j_2,j_3)[X]$. We have $R_p = F_{2,p}(j_{1,p})$ and
$Q_p = F_{3,p}(j_{1,p})$.

We have $P_p = \prod_{M \in S} (X-j_{1,p}(M\tau))$ and with
$$
\widetilde F_{k,p} = \sum_{M \in S} \Biggl( \prod_{B \in S \atop B \not = M}
{X - j_{1,p}(B\tau)} \Biggr) j_{2,p}(M\tau) \in \Q(j_1,j_2,j_3)[X]
$$
we have $R_p = \widetilde F_{2,p}(j_{1,p}) / P_p'(j_{1,p})$ and
$Q_p = \widetilde F_{3,p}(j_{1,p}) / P_p'(j_{1,p})$. Here, $P_p'$ denotes the
derivative of $P_p$. We deduce that it suffices to compute the 3
polynomials $P_p, \widetilde F_{2,p}$ and $\widetilde F_{3,p}$.
In Lemma 6.2 below we prove that the denominators of the coefficients of these
polynomials are closely related to $(p,p)$-{\it split Jacobians\/}.

We say that a principally polarized abelian variety $A/\C$ is $(p,p)$-split
if there exists an isogeny of degree $p^2$ between $A$ and the product
$E \times E'$ of two elliptic curves with the product polarization. The
locus of such $A$ is denoted by ${\Cal L}_p$. It is well known that ${\Cal L}_p$
is a 2-dimensional algebraic subvariety of the 3-dimensional moduli space
$\Am_2$. We have chosen coordinates $j_1,j_2,j_3$ for $\Am_2$, and ${\Cal L}_p$
can be given by an equation $L_p = 0$ for a polynomial
$L_p \in \Q[j_1,j_2,j_3]$.\par
\ \par
\noindent
\proclaim{Lemma 6.2}The denominators of the coefficients of $P_p$,
$\widetilde F_{2,p}$ and $\widetilde F_{3,p}$ are all divisible by the
polynomial $L_p$ describing the moduli space of $(p,p)$-split Jacobians.
\endproclaim
\noindent
{\bf Proof.} Let $\tau\in\H_2$ correspond to a $(p,p)$-split Jacobian, and
let $c$ be a coefficient of~$P_p$.
For some $M \in \Sp(4,\Z) / \Gs2(p)$, the value $j_{1,p}(M\tau)$ is infinite
because the functions $j_i$ have poles at products of elliptic curves.
The evaluation of $c$ at $\tau$ is a symmetric expression in the
$j_{1,p}(M\tau)$'s. Generically, there is no algebraic relation between these
values, and the evaluation of $c$ at $\tau$ is therefore {\it infinite\/}.

Since $j_i(\tau)$ is finite, the
numerator of $c$ is finite. We conclude that the denominator of $c$ must
vanish at $\tau$, i.e., $c$ is divisible by $L_p$. The proof for
$\widetilde F_{2,p}$ and $\widetilde F_{3,p}$ proceeds similarly. \endproof\par
\ \par
\noindent
{\bf Remark.} Points on the variety ${\Cal L_p}$ are in bijection with points
on the {\it Humbert surface\/} $H_{p^2}$, see [\Mura]. It is a traditionally
hard problem to compute equations for Humbert surfaces. Up to now, this has
only been done for $p=2,3,5$, see [\Sha], [\ShaA] and [\ShaB]. As computing
modular polynomials is an even harder
problem, it seems unlikely that we will be able to compute many examples in
the near future. More precisely, we conjecture that the following is true.\par
\ \par\noindent
\proclaim{Conjecture 6.3}There exists an element $c \in \Q_{>0}$ with
the property that the polynomial $P_p$ requires at least $c p^{12}$
bits to write down.
\endproclaim
\noindent
{\bf Heuristic reason.} The degree of $P_p$ equals $(p^4-1)/(p-1) \approx
p^3$ by Lemma~6.1. The numerator of every coefficient of $P_p$ is a polynomial
in $j_1,j_2,j_3$. We expect that the degree in each of the three variables 
$j_1,j_3,j_3$ of these polynomials is of order~$p$. Furthermore, by 
looking at the 1-dimensional case it seems reasonable to assume that these
polynomials themselves have coefficients of 
size at least~$p$. This means that we need at least $cp^4$ bits to store one
coefficient of $P_p$ for some constant $c \in \Q_{>0}$. \endproof\par
\ \par\noindent
In the case $p=2$ it is relatively straightforward to compute the polynomial
$L_2$ describing $(2,2)$-split Jacobians. We refer to [\Sha] for its
construction. We have
\medskip
\halign{\thinspace\thinspace\thinspace $#$\thinspace&#& \thinspace $#$\hfill\cr
L_2&=\thinspace& 236196j_1^5 - 972j_1^4j_2^2 + 5832j_1^4j_2j_3 + 19245600j_1^4j_2 - 8748j_1^4j_3^2\cr
\noalign{\smallskip}
&& - 104976000j_1^4j_3 + 125971200000j_1^4 + j_1^3j_2^4 - 12j_1^3j_2^3j_3 - 77436j_1^3j_2^3\cr
\noalign{\smallskip}
&& + 54j_1^3j_2^2j_3^2 + 870912j_1^3j_2^2j_3 - 507384000j_1^3j_2^2 - 108j_1^3j_2j_3^3 - 3090960j_1^3j_2j_3^2\cr
\noalign{\smallskip}
&& + 2099520000j_1^3j_2j_3 + 81j_1^3j_3^4 + 3499200j_1^3j_3^3 + 78j_1^2j_2^5 - 1332j_1^2j_2^4j_3\cr
\noalign{\smallskip}
&& + 592272j_1^2j_2^4 + 8910j_1^2j_2^3j_3^2 - 4743360j_1^2j_2^3j_3 - 29376j_1^2j_2^2j_3^3 + 9331200j_1^2j_2^2j_3^2\cr
\noalign{\smallskip}
&& + 47952j_1^2j_2j_3^4 - 31104j_1^2j_3^5 - 159j_1j_2^6 + 1728j_1j_2^5j_3 - 41472j_1j_2^5\cr
\noalign{\smallskip}
&& - 6048j_1j_2^4j_3^2 + 6912j_1j_2^3j_3^3 + 80j_2^7 - 384j_2^6j_3.\cr}
\medskip
\noindent
In the remainder in this section we describe the explicit computation of the
entire polynomial $P_2$. Our idea is to use an {\it interpolation\/}
technique, i.e., compute $P_2(j_1(\tau),j_2(\tau),j_3(\tau))\in \C[X]$ for
sufficiently many $\tau\in\H_2$ and
use that information to reconstruct the coefficients of~$P_2$. Unfortunately,
we need to know the full denominators of the coefficients of $P_2$ for this
approach to work. The interpolation problem was also considered by Dupont
in his thesis [\Du].  Without
knowing Lemma~6.2, he succeeded in computing~$P_2$. The approach outlined
below is inspired by Dupont's ideas.

Let $c(j_1,j_2,j_3)$ be a coefficient of $P_2$. The first step is to compute
the degree in $j_1,j_2,j_3$ of its numerator $n(c)$ and its
denominator~$d(c)$. We
{\it fix\/} $y,z\in \Q(i)$ and for a collection of values $x_k \in \Q(i)$ we compute a
value $\tau_k\in\H_2$ with
$$
(j_1(\tau_k),j_2(\tau_k),j_3(\tau_k)) = (x_k,y,z)
$$
by first using Mestre's algorithm [\Mes] to find a genus~2 curve $C$ whose
Igusa invariants are $(x_k,y,z)$ and then finding $\tau_k\in\H_2$
corresponding to $\Jac(C)$. We can therefore evaluate the {\it univariate\/}
rational function $c(x_k,y,z)$. It is now an easy matter to determine the
degree in $j_1$ of $n(c)$ and $d(c)$. Indeed, we check for
which values of $m,n$ the matrix\par
\medskip
$$
M(m,n)=\left(
\vbox{
\halign{$#$\thinspace&\thinspace$#$\thinspace&\thinspace$#$\thinspace&\thinspace\hfil$#$\hfil\thinspace&\thinspace\hfil$#$\hfil\thinspace&\thinspace\hfil$#$\hfil\cr
1&\ldots&x_1^m&-c(x_1,y,z)&\ldots&-c(x_1,y,z)x_1^n\cr
\vdots&\ddots&\vdots&\vdots&\vdots&\vdots\cr
1&\ldots&x_{m+n+2}^m&-c(x_{m+n+2},y,z)&\ldots&-c(x_{m+n+2},y,z)x_{m+n+2}^n\cr}\vskip-.5cm}
\right)
$$
has non-zero solution-space for arbitrary $x_k\in\Q(i)$. The smallest $m,n$
for which this is the case are the degrees of $n(c)$ and~$d(c)$.
We find for instance that the constant term of $P_2$ has a numerator of
degree 60 and a denominator of degree~$51$ in~$j_1$. Likewise, we can find
the degrees in $j_2$ and~$j_3$. The degree in $j_2$ of the denominators is
42 for all coefficients and we find 30 for the degree in~$j_3$.

As $L_2$ has degree $7$ in $j_2$ and degree $5$ in $j_3$, we {\it guess\/} that
$L_2^6$ divides $d(c)$. We are still missing a polynomial in
$j_1$ of degree $>1$ for the denominator. To find this polynomial, a natural
idea is to try $j_1^\alpha$ with $\alpha$ the difference between the degree of $d(c)$
and $6\cdot 5 = 30$. One heuristic reason for this is the
following: if $\tau\in\H_2$ corresponds to the product of elliptic curves, then
the numerator of $c$ is infinite at~$\tau$. Combined with the vanishing of
$L_2^6$ at $\tau$, this would mean that $c$ has a pole of very high order at
such~$\tau$. To `compensate' for this, we multiply $L_2^6$ by ~$j_1^\alpha$.
One can verify that the denominator is indeed $j_1^\alpha L_2^6$ by
taking $x,y,z\in\Z[i]$ and looking at the denominator of $c(x,y,z) \in \Q(i)$.

Having computed the denominator of $c$, it is an easy matter to compute the
numerator. Indeed, we can evaluate $d(c) c$ at any point $\tau\in\H_2$ and
apply interpolation techniques to find~$n(c)$. As the degrees in $j_1$,
$j_2$ and $j_3$ of $n(c)$ are relatively large, this does take a large amount
of time. The constant term of $P_2$ contains $16795$ monomials for instance,
with coefficients up to 200 decimal digits. It takes more than 50 megabytes
to store $P_2, Q_2, R_2$.

Combined with Conjecture~6.3, the fact that $P_2, Q_2, R_2$ require more
than 50 megabytes to store means that the computation of $P_3, Q_3, R_3$ will
be a major effort. We can guess the cost of the computation by combining the
asymptotic growth predicted by Conjecture~6.3 with the `initial data'
provided by $P_2$. A quick computation yields that we expect that every
coefficient of the numerator of $P_2$ has roughly 56000 monomials with
coefficients up to 300 decimal digits. As $P_3$ has degree $(3^4-1)/(3-1)
= 40$, this means that we would roughly need 0.7 gigabytes to store the
numerator of $P_3$. The denominator is fairly small, so we would need
about 2 gigabytes of memory to store $P_3, Q_3, R_3$. This is still
feasible with currect computers, but we did not attempt the computation.

\head 7. Richelot isogeny
\endhead
\noindent
In this Section we zoom in on the special case $p=2$, and explain the link
between our modular polynomials $P_2, Q_2, R_2$ and the classical `Richelot
isogeny'. The Richelot isogeny is typically used as a tool for computing the
Mordell-Weil group of a Jacobian of a genus 2 curve, and we first explain
the construction.

Fix a non-singular curve $C/\C$ of genus $2$, and pick an equation
$$
Y^2 = f(X)
$$
for $C$, where $f \in \C[X]$ is monic of degree~6. We pick a factorization
$f = A B C$ of $f$ into three monic polynomials, each of degree $2$.
Writing $[A,B] = A'B-AB'$, with $A'$ the derivative of $A$, we define the
curve $C'$ by the equation
$$
\Delta Y^2 = [A,B] [A,C] [B,C]. \eqno(7.1)
$$
Here, $\Delta$ is the determinant of $A,B,C$ with respect to the basis
$1,X,X^2$. A simple check shows that the right hand side is again a monic
polynomial of degree $6$. It is separable if and only if $\Delta$ is
non-zero.\par
\ \par
\proclaim{Lemma 7.1} The Jacobians of the curves $C$ and $C'$ defined above
are $(2,2)$-isogenous. Furthermore, every $(2,2)$-isogeny from $\Jac(C)$ into
some principally polarized abelian variety $A$ arises via this construction.
\endproclaim
\noindent
{\bf Proof.} The fact that $\Jac(C)$ and $\Jac(C')$ are $(2,2)$-isogenous
can be found in [\BM]. To show that every $(2,2)$-isogeny is a `Richelot
isogeny', we look at the generic case that $\Jac(C)$ is {\it not\/}
$(2,2)$-split. It then suffices to prove that there are $[\Sp_4(\Z) : \Gs2(2)]
= 15$ different equations for the curve $C'$. This is simple combinatorics:
we have a priori $6\cdot 5 /2 = 15$ choices for the polynomial $A$, and
then $6$ choices for the polynomial~$B$. As the right hand side of (7.1) is
invariant under a permutation of $A,B,C$, we get $15$ different
equations.\endproof\par
\ \par\noindent
The connection with the modular polynomial $P_2$ defined in this paper is as follows.
Assume that $\Jac(C)$ is not $(2,2)$-split. The discriminant $\Delta$ is then
non-zero for every choice of factorization $f = A B C$. We let $\varphi$
be the map sending an Igusa invariant $j_i$ to $j_i(\Jac(C)) \in \C$.
Lemma 7.1 tells us that the $15$ roots of $\varphi(P_2)\in\C[X]$ are exactly the
first Igusa invariants of the curves $C'$ in (7.1).
There are similar relations for $\varphi(R_2)$ and $\varphi(Q_2)$.

\enddocument